\documentclass[12pt]{amsart}
\usepackage[T1]{fontenc}
\usepackage{geometry}
\usepackage{amssymb}

\makeatletter

\providecommand{\tabularnewline}{\\}

 \theoremstyle{plain}    
 \newtheorem{thm}{Theorem}[section]
 \numberwithin{equation}{section} 
 \numberwithin{figure}{section} 
 \theoremstyle{plain}
 \theoremstyle{definition}
 \newtheorem{defn}[thm]{Definition}
 \theoremstyle{plain}    
 \newtheorem{lem}[thm]{Lemma} 
 \theoremstyle{plain}    
 \newtheorem{prop}[thm]{Proposition} 
 \theoremstyle{plain}    
 \newtheorem{cor}[thm]{Corollary} 

\usepackage[all]{xy}

\makeatother
\begin{document}

\title{Rubbling and Optimal Rubbling of Graphs}

\author{Christopher Belford}

\author{N\'Andor Sieben}

\begin{abstract}
A pebbling move on a graph removes two pebbles at a vertex and adds
one pebble at an adjacent vertex. Rubbling is a version of pebbling
where an additional move is allowed. In this new move one pebble is
removed at vertices $v$ and $w$ adjacent to a vertex $u$ and an
extra pebble is added at vertex $u$. A vertex is reachable from a
pebble distribution if it is possible to move a pebble to that vertex
using rubbling moves. The rubbling number of a graph is the smallest
number $m$ needed to guarantee that any vertex is reachable from
any pebble distribution of $m$ pebbles. The optimal rubbling number
is the smallest number $m$ needed to guarantee a pebble distribution
of $m$ pebbles from which any vertex is reachable. We determine the
rubbling and optimal rubbling number of some families of graphs including
cycles.
\end{abstract}

\address{Northern Arizona University, Department of Mathematics and Statistics,
Flagstaff AZ 86011-5717, USA}

\email{cbelford@gmail.com}

\email{nandor.sieben@nau.edu}



\date{\the\month/\the\day/\the\year}

\maketitle

\section{Introduction}

Graph pebbling has its origin in number theory. It is a model for
the transportation of resources. Starting with a pebble distribution
on the vertices of a simple connected graph, a \emph{pebbling move}
removes two pebbles from a vertex and adds one pebble at an adjacent
vertex. We can think of the pebbles as fuel containers. Then the loss
of the pebble during a move is the cost of transportation. A vertex
is called \emph{reachable} if a pebble can be moved to that vertex
using pebbling moves. There are several questions we can ask about
pebbling. How many pebbles will guarantee that every vertex is reachable,
or that all vertices are reachable at the same time? How can we place
the smallest number of pebbles such that every vertex is reachable?
For a comprehensive list of references for the extensive literature
see the survey papers \cite{Hurlbert_survey1,Hurlbert_survey2}. 

In the current paper we propose the study of an extension of pebbling
called \emph{rubbling}. In this version we also allow a move that
removes a pebble from the vertices $v$ and $w$ that are adjacent
to a vertex $u$, and adds a pebble at vertex $u$. We find rubbling
versions of some of the well known pebbling tools such as the transition
digraph, the No Cycle Lemma, squishing and smoothing. We use these
tools to find the rubbling number and the optimal rubbling number
for some families of graphs including complete graphs, complete bipartite
graphs, paths, wheels and cycles.

\section{Preliminaries}

Let $G$ be a simple graph. We use the notation $V(G)$ for the vertex
set and $E(G)$ for the edge set. A \emph{pebble function} on a graph
$G$ is a function $p:V(G)\to{\bf Z}$ where $p(v)$ is the number
of pebbles placed at $v$. A \emph{pebble distribution} is a nonnegative
pebble function. The \emph{size} of a pebble distribution $p$ is
the total number of pebbles $\sum_{v\in V(G)}p(v)$. We are going
to use the notation $p(v_{1},\ldots,v_{n},*)=(a_{1},\ldots,a_{n},q(*))$
to indicate that $p(v_{i})=a_{i}$ for $i\in\{1,\ldots,n\}$ and $p(w)=q(w)$
for all $w\in V(G)\setminus\{ v_{1},\ldots,v_{n}\}$.

\begin{defn}
Consider a pebble function $p$ on the graph $G$. If $\{ v,u\}\in E(G)$
then the \emph{pebbling move} $(v,v\to u)$ removes two pebbles at
vertex $v$ and adds one pebble at vertex $u$ to create a new pebble
function\[
p_{(v,v\to u)}(v,u,*)=(p(v)-2,p(u)+1,p(*)).\]
If $\{ w,u\}\in E(G)$ and $v\not=w$ then the \emph{strict rubbling
move} $(v,w\to u)$ removes one pebble each at vertices $v$ and $w$
and adds one pebble at vertex $u$ to create a new pebble function\[
p_{(v,w\to u)}(v,w,u,*)=(p(v)-1,p(w)-1,p(u)+1,p(*)).\]
A \emph{rubbling move is} either a pebbling move or a strict rubbling
move. 
\end{defn}
Note that the rubbling moves $(v,w\to u)$ and $(w,v\to u)$ are the
same. Also note that the resulting pebble function might not be a
pebble distribution even if $p$ is. 

\begin{defn}
A \emph{rubbling sequence} is a finite sequence $s=(s_{1},\ldots,s_{k})$
of rubbling moves. The pebble function gotten from the pebble function
$p$ after applying the moves in $s$ is denoted by $p_{s}$. 
\end{defn}
The concatenation of the rubbling sequences $r=(r_{1},\ldots,r_{k})$
and $s=(s_{1},\ldots,s_{l})$ is denoted by $rs=(r_{1},\ldots,r_{k},s_{1},\ldots,s_{l})$.

\begin{defn}
A rubbling sequence $s$ is \emph{executable} from the pebble distribution
$p$ if $p_{(s_{1},\ldots,s_{i})}$ is nonnegative for all $i$. A
vertex $v$ of $G$ is \emph{reachable} from the pebble distribution
$p$ if there is an executable rubbling sequence $s$ such that $p_{s}(v)\ge1$.
The \emph{rubbling number} $\rho(G)$ of a graph $G$ is the minimum
number $m$ such that every vertex of $G$ is reachable from any pebble
distribution of size $m$.
\end{defn}
A vertex is reachable if a pebble can be moved to that vertex using
rubbling moves with actual pebbles without ever running out of pebbles.
Changing the order of moves in an executable rubbling sequence $s$
may result in a sequence $r$ that is no longer executable. On the
other hand the ordering of the moves has no effect on the resulting
pebble function, that is, $p_{s}=p_{r}$. This justifies the following
definition.

\begin{defn}
Let $S$ be a multiset of rubbling moves. The pebble function gotten
from the pebble function $p$ after applying the moves in $S$ in
any order is denoted by $p_{S}$.
\end{defn}

\section{The transition digraph and the No Cycle Lemma}

\begin{defn}
Given a multiset $S$ of rubbling moves on $G$, the \emph{transition
digraph} $T(G,S)$ is a directed multigraph whose vertex set is $V(G)$,
and each move $(v,w\to u)$ in $S$ is represented by two directed
edges $(v,u)$ and $(w,u)$. The transition digraph of a rubbling
sequence $s=(s_{1},\ldots,s_{n})$ is $T(G,s)=T(G,S)$, where $S=\{ s_{1},\ldots,s_{n}\}$
is the multiset of moves in $s$. Let $d_{T(G,S)}^{-}$ represent
the in-degree and $d_{T(G,S)}^{+}$ the out-degree in $T(G,S)$. We
simply write $d^{-}$ and $d^{+}$ if the transition digraph is clear
from context.
\end{defn}
The transition digraph only depends on the rubbling moves and the
graph but not on the pebble distribution or on the order of the moves.
It is possible that $T(G,S)=T(G,R)$ even if $S\not=R$. If $T(G,S)=T(G,R)$
then $p_{S}=p_{R}$, so the effect of a rubbling sequence on a pebble
function only depends on the transition digraph. In fact we have the
following. 

\begin{lem}
If $p$ is a pebble function on $G$ and $S$ is a multiset of rubbling
moves then\[
p_{S}(v)=p(v)+d^{-}(v)/2-d^{+}(v)\]
for all $v\in V(G)$.
\end{lem}
\begin{proof}
The three terms on the right hand side represent the original number
of pebbles, the number of pebbles arrived at $v$ and the number of
pebbles moved away from $v$. 
\end{proof}
We are often interested in the value of $q_{R}(v)-p_{S}(v)$. The
function $\Delta$ defined in the following lemma is going to simplify
our notation. The three parameters of $\Delta$ represent the change
in the number of pebbles, the change in the in-degree and the change
in the out-degree. The proof is a trivial calculation.

\begin{lem}
Define $\Delta(a,b,c)=a+b/2-c$. Then\[
q_{R}(v)-p_{S}(v)=\Delta(q(v)-p(v),d_{T(G,R)}^{-}(v)-d_{T(G,S)}^{-}(v),d_{T(G,R)}^{+}(v)-d_{T(G,S)}^{+}(v)).\]
 
\end{lem}
If the rubbling sequence $s$ is executable from a pebble distribution
$p$ then we must have $p_{s}\ge0$. This motivates the following
terminology.

\begin{defn}
A multiset $S$ of rubbling moves on $G$ is \emph{balanced} with
a pebble distribution $p$ \emph{at vertex} $v$ if $p_{S}(v)\ge0$.
We say $S$ is \emph{balanced} with $p$ if $S$ is balanced with
$p$ at all $v\in V(G)$, that is, $p_{S}\ge0$. We say that a rubbling
sequence \emph{}$s$ is balanced with $p$ if the multiset of moves
in $s$ is balanced with $p$.
\end{defn}
$S$ is trivially balanced with a pebble distribution at $v$ if $d_{T(G,S)}^{+}(v)=0$.
The balance condition is necessary but not sufficient for a rubbling
sequence to be executable. The pebble distribution $p(u,v,w)=(1,1,1)$
on the cycle $C_{3}$ is balanced with $s=((u,u\to v),(v,v\to w),(w,w\to u))$,
but $s$ is not executable. The problem is caused by the cycle in
the transition digraph. The goal of this section is to overcome this
difficulty.

\begin{defn}
A multiset of rubbling moves or a rubbling sequence is called \emph{acyclic}
if the corresponding transition digraph has no directed cycles. Let
$S$ be a multiset of rubbling moves. An acyclic multiset $R\subseteq S$
is called an \emph{untangling} of $S$ if $p_{R}\ge p_{S}$.
\end{defn}
\begin{prop}
\label{pro:unfolding}Every multiset of rubbling moves has an untangling.
\end{prop}
\begin{figure}
\begin{center}~\input{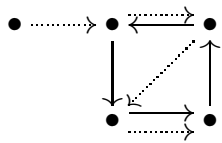}\end{center}

\caption{\label{cap:Arrows-representing-moves} Arrows of $T(G,Q)$. The solid
arrows belong to $C$. \protect \\
}
\end{figure}

\begin{proof}
Let $S$ be the multiset of rubbling moves. Suppose that $T(G,S)$
has a directed cycle $C$. Let $Q$ be the multiset of elements of
$S$ corresponding to the arrows of $C$, see Figure~\ref{cap:Arrows-representing-moves}.
We show that $p_{R}\ge p_{S}$ where $R=S\setminus Q$. If $v\in V(C)$
then there is an $a\le-1$ such that\[
p_{R}(v)-p_{S}(v)=\Delta(0,-2,a)=-1-a\ge0.\]
If $v\in V(G)\setminus V(C)$ then there is an $a\le0$ such that\[
p_{R}(v)-p_{S}(v)=\Delta(0,0,a)\ge0.\]
We can repeat this process on $R$ until we eliminate all the cycles.
This can be finished in finitely many steps since every step decreases
the number of edges in $R$. The resulting multiset is an untangling
of $S$.
\end{proof}
Note that a multiset of moves can have several untanglings. Also note
that if a pebble distribution $p$ is balanced with $S$ and $R$
is an untangling of $S$ then $p_{R}\ge p_{S}\ge0$ and so $p$ is
also balanced with $R$.

\begin{lem}
\label{lem:sourse}If $p$ is a pebble distribution on $G$ that is
balanced with the multiset $S$ of moves and $t=(v,w\to u)\in S$
such that $d^{-}(v)=0=d^{-}(w)$ then $t$ is executable from $p$.
\end{lem}
\begin{proof}
If $v\not=w$ then $p(v)\ge d^{+}(v)\ge1$ and $p(w)\ge d^{+}(w)\ge1$.
If $v=w$ then $p(v)\ge d^{+}(v)\ge2$. In both cases $s$ is executable
from $p$. 
\end{proof}
\begin{prop}
\label{pro:orderability}If the pebble distribution $p$ on $G$ is
balanced with the acyclic multiset $S$ of rubbling moves then there
is a sequence $s$ of the elements of $S$ such that $s$ is executable
from $p$.
\end{prop}
\begin{proof}
We define $s$ recursively. Let $R_{1}=S$. Since $R_{1}$ is acyclic,
we must have a move $s_{1}=(v_{1},w_{1}\to u_{1})\in R_{1}$ such
that $d_{T(G,R_{1})}^{-}(v_{1})=0=d_{T(G,R_{1})}^{-}(w_{1})$. Then
$s_{1}$ is executable from $p$ by Lemma~\ref{lem:sourse}. Let
$R_{i}=R_{i-1}\setminus\{ s_{i-1}\}$. Then $R_{i}$ is acyclic so
we must have a move $s_{i}=(v_{i},w_{i}\to u_{i})\in R_{i}$ such
that $d_{T(G,R_{i})}^{-}(v_{i})=0=d_{T(G,R_{i})}^{-}(w_{i})$. Then
$p_{(s_{1},\ldots,s_{i-1})}$ is balanced with $R_{i}$ since $(p_{(s_{1},\ldots,s_{i-1})})_{R_{i}}=p_{S}\ge0$
and so $s_{i}$ is executable from $p_{(s_{1},\ldots,s_{i-1})}$.
The sequence $s=(s_{1},\ldots,s_{|S|})$ is an ordering of the elements
of $S$ that is executable from $p$.
\end{proof}
The following is the rubbling version of the No-Cycle Lemma for pebbling
\cite{Betsy,Milans,Moews}.

\begin{lem}
\emph{(No Cycle)} Let $p$ be a pebble distribution on $G$ and $v\in V(G)$.
The following are equivalent.
\end{lem}
\begin{enumerate}
\item $v$ is reachable from $p$. 
\item There is a multiset $S$ of rubbling moves such that $S$ is balanced
with $p$ and $p_{S}(v)\ge1$.
\item There is an acyclic multiset $R$ of rubbling moves such that $R$
is balanced with $p$ and $p_{R}(v)\ge1$.
\item $v$ is reachable from $p$ through an acyclic rubbling sequence.
\end{enumerate}
\begin{proof}
If $v$ is reachable from $p$ then there is an executable sequence
$s$ of rubbling moves. The multiset $S$ of rubbling moves of $s$
is balanced with $p$ and $p_{S}(v)\ge1$. So (1) implies (2). If
$S$ satisfies (2) then an untangling $R$ of $S$ satisfies (3).
Suppose $R$ satisfies (3). By Proposition~\ref{pro:orderability},
there is an executable ordering $r$ of the moves of $R$. This $r$
is acyclic and $v$ is reachable through $r$ since $p_{r}(v)=p_{R}(v)\ge1$.
So (3) implies (4). Finally, (4) clearly implies (1).
\end{proof}
\begin{cor}
\label{cor:no-flip-flop}If a vertex is reachable from a pebble distribution
$p$ on $G$ then it is also reachable by a rubbling sequence in which
no move of the form $(v,a\to u)$ is followed by a move of the form
$(u,b\to v)$.
\end{cor}

\section{Basic results}

It is clear from the definition that for all graphs $G$ we have $\rho(G)\le\pi(G)$
where $\pi$ is the pebbling number. For the pebbling number we have
$2^{\text{{\rm diam}}(G)}\le\pi(G)$. This is also true for the rubbling
number. To see this we need to find the rubbling number of a path
first.

\begin{prop}
\label{pro:path}The rubbling number of the path with $n$ vertices
is $\rho(P_{n})=2^{n-1}$.
\end{prop}
\begin{proof}
Let $v_{1},\ldots,v_{n}$ be the consecutive vertices of $P_{n}$.
Let $p(v_{n},*)=(m,0)$ be a pebble distribution from which $v_{1}$
is reachable through the acyclic rubbling sequence $s$. We show that
$m\ge2^{n-1}$. Since $v_{1}$ is reachable and $p(v_{1})=0$, the
balance condition at $v_{1}$ implies that $T(G,s)$ has at least
2 arrows from $v_{2}$ to $v_{1}$ and so $d^{+}(v_{2})\ge2$. Since
$T(G,s)$ has no cycles, there are no arrows from $v_{1}$ to $v_{2}$.
The balance condition at $v_{2}$ now implies that $T(G,s)$ has at
least 4 arrows from $v_{3}$ to $v_{2}$ and so $d^{+}(v_{3})\ge2^{2}$.
An inductive argument shows that $d^{+}(v_{n})\ge2^{n-1}$ and $d^{-}(v_{n})=0$.
The balance condition at $v_{n}$ implies that $m\ge d^{+}(v_{n})\ge2^{n-1}$.
This shows that $2^{n-1}\le\rho(P_{n})$.

It is known \cite{Hurlbert_survey1} that $\pi(P_{n})=2^{n-1}$. The
result now follows from the inequality $2^{n-1}\le\rho(P_{n})\le\pi(P_{n})=2^{n-1}$.
\end{proof}
\begin{figure}
\begin{center}~\input{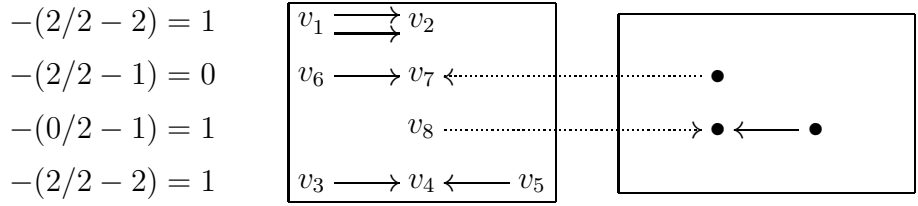}\end{center}

\caption{\label{cap:quotient}Arrows in $T(G,S)$ representing the possible
types of rubbling moves in $E$. The vertices in the same box are
equivalent. The solid arrows connect equivalent vertices. The calculation
on the left shows the change in $\sum_{i}(\frac{1}{2}d^{-}(v_{i})-d^{+}(v_{i}))$
after the removal of one of the rubbling moves. }
\end{figure}

\begin{prop}
If the graph $G$ has diameter $d$ then $2^{d}\le\rho(G)$.
\end{prop}
\begin{proof}
Let $v_{0}$ and $v_{d}$ be vertices at distance $d$. Let $p(v_{0},*)=(m,0)$
be a pebble distribution from which $v_{d}$ is reachable through
the rubbling sequence $s$. We now build a quotient rubbling problem.
Let $[v]$ be the equivalence class of $v$ in the partition of the
vertices of $G$ according to their distances from $v_{0}$. The quotient
simple graph $H$ is isomorphic to $P_{d+1}$ with leafs $[v_{0}]=\{ v_{0}\}$
and $[v_{d}]$. Let $q([v])=\sum_{w\in[v]}p(w)$ for all $[v]\in V(H)$
and note that $q([v_{0}],*)=(m,0)$. The rubbling sequence $s$ induces
a multiset $R$ of rubbling moves on $H$. We construct this $R$
from the multiset $S$ of rubbling moves of $s$. Let $E$ be the
multiset of moves of $S$ of the form $(v,w\to u)$ where $v\in[u]$
or $w\in[u]$. Define $R$ to be the multiset of moves of the form
$([v],[w]\to[u])$ where $(v,w\to u)$ runs through the elements of
$S\setminus E$. 

We show that $R$ is balanced with $q$ . Figure~\ref{cap:quotient}
shows the possible types of moves in $E$. The removal of any of these
moves does not decrease the value of $\sum_{v_{i}\in[v]}(\frac{1}{2}d^{-}(v_{i})-d^{+}(v_{i}))$
and so\[
q_{R}([v])=\sum_{v_{i}\in[v]}p_{S\setminus E}(v_{i})\ge\sum_{v_{i}\in[v]}p_{S}(v_{i})\ge0\]
since $p$ is balanced with $S$. 

We also have $q_{R}([v_{d}])\ge1$ since $v_{d}$ is reachable and
so $p_{S}(v_{d})\ge1$. Thus $[v_{d}]$ is reachable from $q$ and
so the result now follows from Proposition~\ref{pro:path}. 
\end{proof}
For the pebbling number we have $\pi(G)\ge|V(G)|$. This inequality
does not hold for the rubbling number as we can see in the next result.

\begin{prop}
We have the following values for the rubbling number:

\emph{a.} $\rho(K_{n})=2$ for $n\ge2$ where $K_{n}$ is the complete
graph with $n$ vertices\emph{;}

\emph{b.} $\rho(W_{n})=4$ for $n\ge4$ \emph{}where $W_{n}$ is the
wheel with $n$ spikes\emph{;}

\emph{c.} $\rho(K_{m,n})=4$ for $m,n\ge2$ \emph{}where $K_{m,n}$
is a complete bipartite graph\emph{;}

\emph{d.} $\rho(Q^{n})=2^{n}$ \emph{}for $n\ge1$ \emph{}where $Q^{n}$
is the $n$-dimensional hypercube\emph{;}

\emph{e.} $\rho(G)=2^{s+1}$ where $s$ is the number of vertices
in the spine of the caterpillar $G$.
\end{prop}
\begin{proof}
a. A single pebble is clearly not sufficient but any vertex is reachable
with two pebbles using a single move.

b. If we have 4 pebbles then we can move 2 pebbles to the center using
two moves. Then any other vertex is reachable from the center in a
single move. On the other hand $\rho(W_{n})\ge2^{\text{diam}(W_{n})}=2^{2}=4$. 

c. It is easy to see that from any pebble distribution of size 4 any
vertex is reachable in at most 3 moves. On the other hand we have
$\rho(K_{m,n})\ge2^{\text{diam}(K_{m,n})}=2^{2}=4$. 

d. We know \cite{Chung} that $\pi(Q^{n})=2^{n}$. The result now
follows from the inequality $2^{n}=2^{\text{diam}(Q^{n})}\le\rho(Q^{n})\le\pi(Q^{n})=2^{n}$.

e. The result follows easily from Proposition~\ref{pro:path}.
\end{proof}
\begin{figure}
\begin{center}~\input{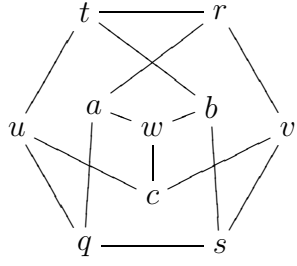}\end{center}

\caption{\label{cap:The-Petersen-graph}The Petersen graph $P$.\protect \\
}
\end{figure}

\begin{prop}
The rubbling number of the Petersen graph $P$ is $\rho(P)=5$.
\end{prop}
\begin{proof}
Consider Figure~\ref{cap:The-Petersen-graph}. It is easy to see
that vertex $w$ is not reachable from the pebble distribution $p(r,s,*)=(3,1,0)$
and so $\rho(P)>4$. To show that $\rho(P)\le5$, assume that a vertex
is not reachable from a pebble distribution $p$ of size 5. Since
$P$ is vertex transitive, we can assume that this vertex is $w$.
Then we must have \[
p(a)+p(b)+p(c)+\left\lfloor \frac{p(q)+p(r)}{2}\right\rfloor +\left\lfloor \frac{p(s)+p(t)}{2}\right\rfloor +\left\lfloor \frac{p(u)+p(v)}{2}\right\rfloor \le1,\]
otherwise we could make the total number of pebbles at vertices $a$,
$b$ and $c$ more than 2 after which $w$ is reachable. This inequality
forces $p(a)=p(b)=p(c)=0$ and two of the remaining terms to be 0
 as well. So by symmetry we can assume that the last term is 1 and
all the other terms are 0. Then we must have $p(u)+p(v)=3$ and $p(q)+p(r)=1=p(s)+p(t)$.
A simple case analysis shows that $w$ is reachable from this $p$,
which is a contradiction. 
\end{proof}

\section{Squishing}

The following terms are needed for the rubbling version of the squishing
lemma of \cite{Bunde_optimal}. A \emph{thread} in a graph is a path
containing vertices of degree 2. A pebble distribution is \emph{squished}
on a thread $P$ if all the pebbles on $P$ are placed on a single
vertex of $P$ or on two adjacent vertices of $P$.

\begin{lem}
\label{lem:notype2}Let $P$ be a thread in $G$. If vertex $x\not\in V(P)$
is reachable from the pebble distribution $p$ then $x$ is reachable
from $p$ through a rubbling sequence in which there is no strict
rubbling move of the form $(v,w\to u)$ where $u\in V(P)$.
\end{lem}
\begin{proof}
Let $S$ be an acyclic multiset of rubbling moves balanced with $p$
such that $p_{S}(x)\ge1$. Let $E$ be the multiset of strict rubbling
moves of $S$ of the form $(v,w\to u)$ where $u\in V(P)$. 

If $e=(v,w\to u)\in E$ then we have $d_{T(G,S\setminus\{ e\})}^{+}(u)=d_{T(G,S)}^{+}(u)=0$
since $S$ is acyclic and so $S\setminus\{ e\}$ is balanced with
$p$ at $u$. It is clear that $p_{S\setminus\{ e\}}(y)\ge p_{S}(y)$
for all $y\in V(G)\setminus\{ u\}$ and so $S\setminus\{ e\}$ is
balanced with $p$. We still know that $S\setminus\{ e\}$ is acyclic
and $p_{S\setminus\{ e\}}(x)\ge1$, so induction shows that $R=S\setminus E$
is balanced with $p$.

By Proposition~\ref{pro:orderability}, there is an ordering $r$
of the elements of $R$ that is executable from $p$. Then $v$ is
reachable through $r$ since $p_{r}(v)=p_{S}(v)\ge1$. 
\end{proof}
The following is the rubbling version of the Squishing Lemma for pebbling
\cite{Bunde_optimal}.

\begin{lem}
\emph{(Squishing)} If vertex $v$ is not reachable from a pebble distribution
with size $n$ then there is a pebble distribution $r$ of size $n$
that is squished on each thread not containing $v$ such that $v$
is not reachable from $r$ either.
\end{lem}
\begin{proof}
The result follows from \cite[Lemma 4]{Bunde_optimal} and \ref{lem:notype2}.
\end{proof}

\section{Rubbling $C_{n}$}

The Squishing Lemma allows us to find the rubbling numbers of cycles.
For the pebbling numbers of $C_{n}$ see \cite{Pachter,Bunde_optimal}.

\begin{prop}
The rubbling number of an even cycle is $\rho(C_{2k})=2^{k}$.
\end{prop}
\begin{proof}
It is well known \cite{Pachter} that $\pi(C_{2k})=2^{k}$. The first
result now follows since\[
2^{k}=2^{\text{diam}(C_{2k})}\le\rho(C_{2k})\le\pi(C_{2k})=2^{k}.\]

\end{proof}
\begin{prop}
The rubbling number of an odd cycle is $\rho(C_{2k+1})=\lfloor\frac{7\cdot2^{k-1}-2}{3}\rfloor+1$.
\end{prop}
\begin{proof}
Let $C_{2k+1}$ be the cycle with consecutive vertices \[
x_{k},x_{k-1},\ldots,x_{1},v,y_{1},y_{2},\ldots,y_{k},x_{k}.\]
 First we show that $\rho(C_{2k+1})\le\lfloor\frac{7\cdot2^{k-1}-2}{3}\rfloor+1$.
Let $p$ be a pebble distribution on $C_{2k+1}$ from which not every
vertex is reachable. It suffices to show that $p$ contains at most
$\lfloor\frac{7\cdot2^{k-1}-2}{3}\rfloor$ pebbles. By symmetry, we
can assume that $v$ is the vertex that is not reachable from $p$.
By the Squishing Lemma, we can assume that $p$ is squished on the
thread with consecutive vertices $y_{1},\ldots,y_{k},x_{k},\ldots,x_{1}$. 

First we consider the case when all the pebbles are at distance $k$
from $v$, that is, $p(x_{k},y_{k},*)=(a,b,0)$. By symmetry, we can
assume that $0\le a\le b$. Then we must have\begin{equation}
\left\lfloor \frac{a}{2}\right\rfloor +b\le2^{k}-1,\label{eq:1}\end{equation}
otherwise we could move $\lfloor\frac{a}{2}\rfloor$ pebbles from
vertex $x_{k}$ to vertex $y_{k}$ and then reach $v$ from $b_{k}$.
Hence $\frac{a}{2}<\left\lfloor \frac{a}{2}\right\rfloor +1\le2^{k}-1-b+1=2^{k}-b$
and so\begin{equation}
a+2b\le2^{k+1}-1.\label{eq:1a}\end{equation}
 We also must have \begin{equation}
\left\lfloor \frac{b-2^{k-1}}{2}\right\rfloor +a\le2^{k-1}-1,\label{eq:2}\end{equation}
otherwise we could move $\lfloor\frac{b-2^{k-1}}{2}\rfloor$ pebbles
from vertex $y_{k}$ to vertex $x_{k}$ after which $x_{1}$ is reachable
from $x_{k}$ and $y_{1}$ is reachable from $y_{k}$, and so $v$
would be reachable by the move $(x_{1},y_{1}\to v)$. Hence $\frac{b-2^{k-1}}{2}<\left\lfloor \frac{b-2^{k-1}}{2}\right\rfloor +1\le2^{k-1}-1-a+1=2^{k-1}-a$
and so \begin{equation}
b+2a\le2^{k}+2^{k-1}-1.\label{eq:2a}\end{equation}
 Adding (\ref{eq:1a}) and (\ref{eq:2a}) gives\[
3(a+b)\le2^{k+1}-1+2^{k}+2^{k-1}-1=7\cdot2^{k-1}-2,\]
which shows that $|p|=a+b\le\lfloor\frac{7\cdot2^{k-1}-2}{3}\rfloor$.

Now we consider the case when some pebbles are closer to $v$ than
$k$, that is, $p(x_{i},x_{i+1},*)=(b,a,0)$ with $b\ge1$ and $a\ge0$
for some $1\le i<k$. Then we must have $\left\lfloor \frac{a}{2}\right\rfloor +b\le2^{i}-1\le2^{k-1}-1$
otherwise $v$ is reachable. Hence\begin{eqnarray*}
|p| & = & a+b\le a-\left\lfloor \frac{a}{2}\right\rfloor +\left\lfloor \frac{a}{2}\right\rfloor +b\\
 & \le & \left\lfloor \frac{a}{2}\right\rfloor +1+2^{k-1}-1\le2^{k-1}-1-b+1+2^{k-1}-1\\
 & = & 2\cdot2^{k-1}-2<\left\lfloor \frac{7\cdot2^{k-1}-2}{3}\right\rfloor .\end{eqnarray*}

Now we show that we can always distribute $\lfloor\frac{7\cdot2^{k-1}-2}{3}\rfloor$
pebbles so that $v$ is unreachable and so $\rho(C_{2k+1})\ge\lfloor\frac{7\cdot2^{k-1}-2}{3}\rfloor+1$.
Let $a=\lfloor\frac{2^{k}}{3}\rfloor$ and $b=\lfloor\frac{5\cdot2^{k-1}}{3}\rfloor$.
It is easy to check that \[
a=\begin{cases}
\frac{2^{k}-2}{3}, & \text{$k$ odd}\\
\frac{2^{k}-1}{3}, & \text{$k$ even}\end{cases},\  b=\begin{cases}
\frac{5\cdot2^{k-1}-2}{3}, & \text{$k$ odd}\\
\frac{5\cdot2^{k-1}-1}{3}, & \text{$k$ even}\end{cases},\ \left\lfloor \frac{7\cdot2^{k-1}-2}{3}\right\rfloor =\begin{cases}
\frac{7\cdot2^{k-1}-4}{3}, & \text{$k$ odd}\\
\frac{7\cdot2^{k-1}-2}{3}, & \text{$k$ even}\end{cases}\]
and so $a+b=\lfloor\frac{7\cdot2^{k-1}-2}{3}\rfloor$. We show that
$v$ is unreachable from the pebble distribution $p(x_{k},y_{k},*)=(a,b,0)$. 

It is easy to see that $a$ and $b$ satisfy (\ref{eq:1a}) and (\ref{eq:2a}).
Suppose that $v$ is reachable from $p$, that is, there is an acyclic
multiset $S$ of rubbling moves that is balanced with $p$ satisfying
$p_{S}(v)\ge1$. The balance condition at $v$ shows that $d^{-}(v)\ge2$.
Hence $S$ must have at least one of $(x_{1},y_{1}\to v)$, $(x_{1},x_{1}\to v)$
or $(y_{1},y_{2}\to v)$. 

First assume that $(x_{1},y_{1}\to v)\in S$. The argument used in
the proof of Proposition~\ref{pro:path} shows that then $T(G,S)$
has at least $2^{i-1}$ arrows from $x_{i}$ to $x_{i-1}$ and from
$y_{i}$ to $y_{i-1}$ for all $i\in\{2,\ldots,k\}$. Since $S$ is
acyclic, any arrow in $T(G,S)$ pointing to $x_{k}$ must come from
$y_{k}$. So the balance condition at $x_{k}$ requires $m$ arrows
from $y_{k}$ to $x_{k}$ satisfying $2^{k-1}\le a+\frac{m}{2}$.
The balance condition at $y_{k}$ gives $2^{k-1}+m\le b$. Combining
the two inequalities gives $2^{k}+2^{k-1}\le b+2a$ which contradicts
(\ref{eq:2a}).

Next assume that $(y_{1},y_{1}\to v)\in S$. Then $T(G,S)$ has at
least $2^{i}$ arrows from $y_{i}$ to $y_{i-1}$ for all $i\in\{2,\ldots,k\}$.
The balance condition at $y_{k}$ requires $m$ arrows from $x_{k}$
to $y_{k}$ satisfying $2^{k}\le b+\frac{m}{2}$. We must have $d^{-}(x_{k})=0$,
otherwise there is a directed path from $v$ to $x_{k}$ which is
impossible since $S$ is acyclic. The balance condition at $x_{k}$
gives $m\le a$. Combining the two inequalities gives $2^{k+1}\le a+2b$
which contradicts (\ref{eq:1a}).

Similar argument shows that $(x_{1},x_{1}\to v)\in S$ is also impossible.
\end{proof}

\section{Optimal rubbling}

Optimal pebbling was studied in \cite{Pachter,Moews_optimal,Fu,Bunde_optimal}.
In this section we investigate the optimal rubbling number of certain
graphs. 

\begin{defn}
The \emph{optimal rubbling number} $\rho_{\text{opt}}(G)$ of a graph
$G$ is the minimum number $m$ for which there is a pebble distribution
of size $m$ from which every vertex of $G$ is reachable.
\end{defn}
\begin{prop}
We have the following values for the optimal rubbling number:

\emph{a.} $\rho_{\text{{\rm opt}}}(K_{n})=2$ for $n\ge2$ where $K_{n}$
is the complete graph with $n$ vertices\emph{;}

\emph{b.} $\rho_{{\rm opt}}(W_{n})=2$ for $n\ge4$ \emph{}where $W_{n}$
is the wheel with $n$ spikes\emph{;}

\emph{c.} $\rho_{\text{{\rm opt}}}(K_{m,n})=3$ for $m,n\ge3$ \emph{}where
$K_{m,n}$ is the complete bipartite graph\emph{;}

\emph{d}. $\rho_{{\rm opt}}(P)=4$ where $P$ is the Petersen graph.
\end{prop}
\begin{proof}
a. Not every vertex of $K_{n}$ is reachable from a distribution of
size 1 since $n\ge2$. On the other hand any vertex is reachable by
a single move from any distribution of size 2.

b. Again, not every vertex of $W_{n}$ is reachable from a distribution
of size 1. On the other hand, every vertex is reachable from the distribution
that has 2 pebbles at the center of $W_{n}$. 

c. Let $A$ and $B$ be the natural partition of the vertex set of
$K_{m,n}$. Let $p$ be a pebble distribution of size 2. If $p$ places
both pebbles on vertices in $A$ then there is a vertex in $A$ that
is not reachable from $p$. If $p$ places both pebbles on vertices
in $B$ then there is a vertex in $B$ that is not reachable from
$p$. If $p$ places one pebble on a vertex in $A$ and one pebble
on a vertex in $B$ then both $A$ and $B$ have vertices that are
unreachable from $p$. On the other hand any vertex is reachable in
at most two moves from a pebble distribution that places one pebble
on a vertex in $A$ and two pebbles on a vertex in $B$.

d. Every vertex is reachable from the pebble distribution that has
4 pebbles on any of the vertices. A simple case analysis shows that
3 pebbles are not sufficient to make every vertex reachable.
\end{proof}
Rolling moves serve the same purpose as the smoothing move of \cite{Bunde_optimal}.

\begin{figure}
\begin{center}~\input{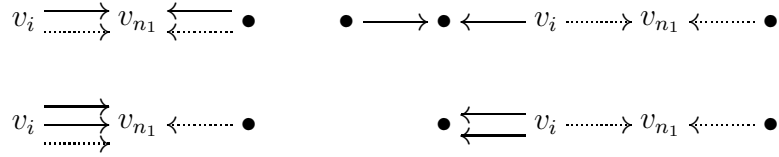}\end{center}

\caption{\label{cap:droll}The four possible configurations for $T(G,S\setminus\tilde{R})$.
The solid arrows represent the moves corresponding to the arrows of
$\tilde{P}$. The dotted arrows represent the moves corresponding
to the arrows of $P$.}
\end{figure}

\begin{defn}
Let $v_{1},\ldots,v_{n}$ be the consecutive vertices of a path such
that the degree of $v_{1}$ is 1 and the degrees of $v_{2},v_{3},\ldots,v_{n-1}$
are all 2. The subgraph induced by $\{ v_{1},\ldots,v_{n}\}$ is called
an \emph{arm} of the graph. Let $p$ be a pebble distribution such
that $p(v_{i})\ge2$ for some $i\in\{1,\ldots,n-1\}$, $p(v_{n})=0$,
and $p(v_{j})\ge1$ for all $j\in\{1,\ldots,n-1\}$. A \emph{single
rolling move} creates a new pebble distribution $q$ by taking one
pebble from $v_{i}$ and placing it on $v_{n}$, that is $q(v_{i},v_{n},*)=(p(v_{i})-1,1,p(*))$.
See Figure~\ref{cap:rollvis}.

\begin{figure}
~\input{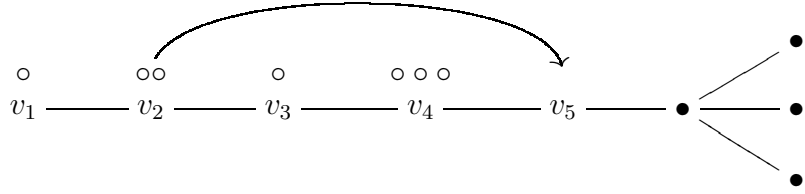}

\caption{\label{cap:rollvis}Visualization of a single rolling move with $i=2$
and $n=5$. An arrow indicates the transfer of a single pebble}
\end{figure}

\end{defn}
\begin{lem}
\label{lem:roll}Let $q$ be a pebble distribution on $G$ gotten
from the pebble distribution $p$ by applying a single rolling move
from $v_{i}$ to $v_{n}$ on the arm with vertices $v_{1},\ldots,v_{n}$.
If vertex $u\in G$ is reachable from $p$ then $u$ is also reachable
from $q$.

\begin{figure}
\begin{center}~\input{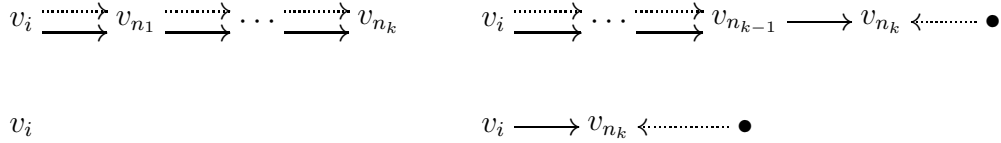}\end{center}

\caption{\label{cap:roll}Four possible configurations for $T(G,S\setminus R)$.
The solid arrows represent the arrows of $P$.}
\end{figure}

\end{lem}
\begin{proof}
If $u$ is a vertex of the arm then it is clearly reachable from $q$
so we can assume that $u$ is not on the arm. Let $S$ be an acyclic
multiset of rubbling moves balanced with $p$ such that $p_{S}(u)\ge1$.
Let $P$ be a maximum length directed path in $T(G,S)$ starting at
$v_{i}$ and not going further than $v_{n}$. Then $P$ has consecutive
vertices $v_{i}=v_{n_{0}},v_{n_{1}}\ldots,v_{n_{k}}$ on the arm.
Let $R$ be the multiset containing the elements of $S$ without the
moves corresponding to the arrows of $P$. We show that $R$ is balanced
with $q$ and so $u$ is reachable from $q$ since $q_{R}(u)=p_{S}(u)\ge1$.
Figure ~\ref{cap:roll} shows the possible configurations for $T(G,S\setminus R)$.
We have $d_{T(G,S)}^{+}(v_{n_{k}})=0$ even if $n_{k}=1$. If $n_{k}=n$
then \[
q_{R}(v_{n_{k}})=p_{S}(v_{n_{k}})+\Delta(1,-2,0)=p_{S}(v_{n_{k}})\ge1\ge0,\]
while if $n_{k}\not=n$ then\[
q_{R}(v_{n_{k}})=p_{S}(v_{n_{k}})+\Delta(0,-2,0)\ge p_{S}(v_{n_{k}})-1\ge2-1\ge0.\]
So $R$ is balanced with $q$ at $v_{n_{k}}$. If $d_{T(G,S)}^{+}(v_{n_{0}})=0$
then $n_{0}=n_{k}$, otherwise there is an $a\in\{-1,-2\}$ such that\[
q_{R}(v_{n_{0}})=p_{S}(v_{n_{0}})+\Delta(-1,0,a)\ge p_{S}(v_{n_{0}})\ge0\]
and so $R$ is balanced with $q$ at $v_{n_{0}}$. If $0<j<k$ then
there is an $a\in\{-1,-2\}$ such that\[
q_{R}(v_{n_{j}})=p_{S}(v_{n_{j}})+\Delta(0,-2,a)\ge p_{S}(v_{n_{j}})\ge0\]
and so $R$ is balanced with $q$ at $v_{n_{j}}.$ It is clear that
$R$ is balanced with $q$ at every other vertex. 
\end{proof}
\begin{defn}
Let $v_{1},\ldots,v_{n}$ be the consecutive vertices of a path such
that the degrees of $v_{2},v_{3},\ldots,v_{n-1}$ are all 2. Let $p$
be a pebble distribution such that $p(v_{1})=0=p(v_{n})$, $p(v_{i})\ge2$
for some $i\in\{2,\ldots,n-1\}$ and $p(v_{j})\ge1$ for all $j\in\{2,\ldots,n-1\}$.
A \emph{double rolling move} creates a new pebble distribution $q$
by taking two pebbles from $v_{i}$ and placing one pebble on $v_{1}$
and one pebble on $v_{n}$, that is $q(v_{i},v_{1},v_{n},*)=(p(v_{i})-2,1,1,p(*))$.
See Figure~\ref{cap:drollvis}.

\begin{figure}
~\input{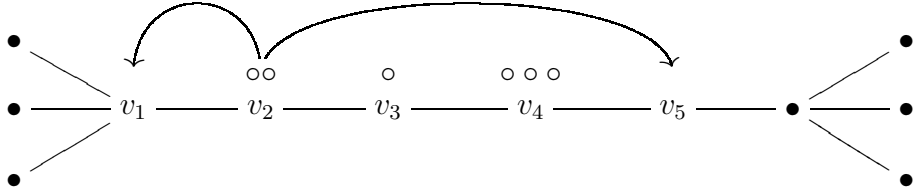}

\caption{\label{cap:drollvis}Visualization of a double rolling move with
$i=2$ and $n=5$. An arrow indicates the transfer of a single pebble.}
\end{figure}

\end{defn}
\begin{lem}
\label{lem:droll}Let $q$ be a pebble distribution on $G$ gotten
from the pebble distribution $p$ by applying a double rolling move
from vertex $v_{i}$ to vertices $v_{1}$ and $v_{n}$ on the path
with consecutive vertices $v_{1},\ldots,v_{n}$. If vertex $u\in G$
is reachable from $p$ then $u$ is also reachable from $q$. 
\end{lem}
\begin{proof}
If $u\in\{ v_{1},\ldots,v_{n}\}$ then it is clearly reachable from
$q$ so we can assume that $u\not\in\{ v_{1},\ldots,v_{n}\}$. Let
$S$ be an acyclic multiset of rubbling moves balanced with $p$ such
that $p_{S}(u)\ge1$. Let $P$ be a maximum length directed path in
$T(G,S)$ starting at $v_{i}$ and not going further than $v_{1}$
or $v_{n}$. Then $P$ has consecutive vertices $v_{i}=v_{n_{0}},v_{n_{1}}\ldots,v_{n_{k}}\in\{ v_{1},\ldots,v_{n}\}$.
Let $R$ be the multiset containing the elements of $S$ without the
moves corresponding to the arrows of $P$. An argument similar to
the one in the proof of Lemma~\ref{lem:roll} shows that $R$ is
clearly balanced with $q$ at every vertex except maybe at $v_{i}$.
If $n_{k}=n_{0}$ or the arrow $(v_{n_{0}},v_{n_{1}})$ in $P$ corresponds
to a pebbling move, then $R$ is balanced with $q$ at $v_{i}$ as
well. Then $u$ is reachable from $q$ since $q_{R}(u)=p_{S}(u)\ge1$. 

So we can assume that $(v_{n_{0}},v_{n_{1}})$ corresponds to a strict
rubbling move and that $k=1$. Let $\tilde{P}$ be a maximum length
path in $T(G,R)$. Since $k=1$, the length of $\tilde{P}$ is either
0 or 1. If this length is 0, then $q$ is balanced with $R$ at $v_{i}$
since $d_{T(G,R)}^{+}(v_{i})=0$ and we are done. If the length of
$\tilde{P}$ is 1, then let $\tilde{R}$ be the multiset containing
the elements of $R$ without the moves corresponding to the arrows
of $\tilde{P}$. Figure~\ref{cap:droll} shows the possibilities
for $T(G,S\setminus\tilde{R})$. It is easy to check that $\tilde{R}$
is balanced with $q$ in each case. Thus $u$ is reachable from $q$
since $q_{\tilde{R}}(u)\ge p_{S}(u)$.
\end{proof}
 Rolling moves make it possible to find the optimal rubbling number
of paths and cycles.

\begin{prop}
\emph{The optimal rubbling number of the path is} $\rho_{\text{{\rm opt}}}(P_{n})=\lceil\frac{n+1}{2}\rceil$.
\end{prop}
\begin{proof}
Let $P_{n}$ be the path with consecutive vertices $v_{1},\ldots,v_{n}$.
It is clear that every vertex is reachable from the pebble distribution\[
p(v_{i})=\begin{cases}
1, & \text{$i$ is odd or $i=n$}\\
0, & \text{else}\end{cases}\]
which has size $\lceil\frac{n+1}{2}\rceil$. 

Now assume that there is a pebble distribution of size $\lceil\frac{n+1}{2}\rceil-1$
from which every vertex of $P_{n}$ is reachable. Let us apply all
available rolling moves (single or double). The process ends in finitely
many steps since a rolling move reduces the number of pebbles on vertices
with more than one pebble by at least one. If there is a vertex with
more than one pebble and a vertex with no pebbles, then a rolling
move is available. The number of pebbles is not larger than the number
of vertices, so the resulting pebble distribution $q$ has at most
one pebble on each vertex. Every vertex of $P_{n}$ still must be
reachable from $q$ by Lemma~\ref{lem:droll}.

The only moves executable directly from $q$ are strict rubbling moves.
By the No Cycle Lemma we can assume that every vertex is reachable
by a sequence of moves in which a strict rubbling move $(x,y\to z)$
is not followed by a move of the form $(z,z\to x)$ or $(z,z\to y)$.
So we can assume that every vertex is reachable through strict rubbling
moves. Then we must have $q(v_{1})=1=q(v_{n})$ otherwise $v_{1}$
or $v_{n}$ is not reachable. A pigeon hole argument shows that there
must be two neighbor vertices $u$ and $w$ such that $q(u)=0=q(w)$.
But then neither $u$ nor $w$ is reachable from $q$, which is a
contradiction.
\end{proof}
\begin{prop}
The optimal rubbling number of the cycle is $\rho_{\text{{\rm opt}}}(C_{n})=\lceil\frac{n}{2}\rceil$
for $n\ge3$.
\end{prop}
\begin{proof}
Let $C_{n}$ be the cycle with consecutive vertices $v_{1},\ldots,v_{n}$.
It is clear that every vertex is reachable from the pebble distribution\[
p(v_{i})=\begin{cases}
1, & \text{$i$ is odd}\\
0, & \text{else}\end{cases}\]
which has size $\lceil\frac{n}{2}\rceil$. 

Now assume that there is a pebble distribution of size $\lceil\frac{n}{2}\rceil-1$
from which every vertex of $C_{n}$ is reachable. Let us apply all
available double rolling moves. The process ends in finitely many
steps since a double rolling move reduces the number of pebbles on
vertices with more than one pebble by two . If there is a vertex with
more than one pebble and two vertices with no pebbles, then a double
rolling move is available. The number of pebbles is smaller than the
number of vertices, so the resulting pebble distribution $q$ has
at most one pebble on each vertex. Every vertex of $C_{n}$ still
must be reachable from $q$.

The only moves executable directly from $q$ are strict rubbling moves.
The No Cycle Lemma implies that we can assume that every vertex is
reachable through strict rubbling moves. A pigeon hole argument shows
that there must be two neighbor vertices $u$ and $w$ such that $q(u)=0=q(w)$.
But then neither $u$ nor $w$ is reachable from $q$ which is a contradiction.
\end{proof}

\section{Further questions}

There are plenty of unanswered questions. The following might not
be too hard to answer. 

\begin{itemize}
\item What is the optimal rubbling number for the hypercube $Q^{n}$. \emph{}It
is fairly easy to get answers for small $n$ with a computer. The
known values are listed in Table~\ref{cap:Known-rho-opt-hyper}.%
\begin{table}
\begin{tabular}{|c|c|c|c|c|}
\hline 
$n$&
2
&
3
&
4
&
5
\tabularnewline
\hline 
$\rho(B_{n})$&
4&
16&
$>23$&
\tabularnewline
\hline 
$\rho_{{\rm opt}}(B_{n})$&
2&
4&
6&
\tabularnewline
\hline 
$\rho_{{\rm opt}}(Q^{n})$&
2
&
3
&
4
&
6
\tabularnewline
\hline
\end{tabular}

~

\caption{\label{cap:Known-rho-opt-hyper}Rubbling values without a known general
formula.\protect \\
}
\end{table}

\item Does Graham's conjecture hold for the rubbling number?
\item Is the cover rubbling number the same as the cover pebbling number
for every graph?
\item We have $\pi(P_{n})=\rho(P_{n})$, $\pi(Q^{n})=\rho(Q^{n})$ and it
is easy to check that $\pi(L)=8=\rho(L)$ where $L$ is the Lemke
graph \cite{Hurlbert_survey2}. This is not always the case though.
Is it possible to characterize those graphs for which the pebbling
and the rubbling numbers are the same?
\item Let $f(d,n)=\max\{\rho(G)\mid|V(G)|=n\text{ and diam}(G)=d\}$. It
is not hard to check that $f(2,n)\le5$ and $f(3,n)\le9$ for $n\in\{1,\ldots,7\}$.
Do these upper limits hold for all $n$? Is it true that $f(d,n)\le2^{d}+1$
for all $d$ and $n$? 
\end{itemize}
\bibliographystyle{amsplain}
\bibliography{pebble}

\providecommand{\bysame}{\leavevmode\hbox to3em{\hrulefill}\thinspace}
\providecommand{\MR}{\relax\ifhmode\unskip\space\fi MR }
\providecommand{\MRhref}[2]{%
  \href{http://www.ams.org/mathscinet-getitem?mr=#1}{#2}
}
\providecommand{\href}[2]{#2}
\begin{thebibliography}{10}

\bibitem{Bunde_optimal}
David~P. Bunde, Erin~W. Chambers, Daniel Cranston, Kevin Milans, and Douglas~B.
  West, \emph{Pebbling and optimal pebbling in graphs}, (Preprint).

\bibitem{Chung}
Fan R.~K. Chung, \emph{Pebbling in hypercubes}, SIAM J. Discrete Math.
  \textbf{2} (1989), no.~4, 467--472.

\bibitem{Betsy}
Betsy Crull, Tammy Cundiff, Paul Feltman, Glenn~H. Hurlbert, Lara Pudwell,
  Zsuzsanna Szaniszlo, and Zsolt Tuza, \emph{The cover pebbling number of
  graphs}, Discrete Math. \textbf{296} (2005), no.~1, 15--23.

\bibitem{Fu}
Hung-Lin Fu and Chin-Lin Shiue, \emph{The optimal pebbling number of the
  complete {$m$}-ary tree}, Discrete Math. \textbf{222} (2000), no.~1-3,
  89--100.

\bibitem{Hurlbert_survey1}
Glenn Hurlbert, \emph{A survey of graph pebbling}, Proceedings of the Thirtieth
  Southeastern International Conference on Combinatorics, Graph Theory, and
  Computing (Boca Raton, FL, 1999), vol. 139, 1999, pp.~41--64.

\bibitem{Hurlbert_survey2}
\bysame, \emph{Recent progress in graph pebbling}, Graph Theory Notes N. Y.
  \textbf{49} (2005), 25--37.

\bibitem{Milans}
Kevin Milans and Bryan Clark, \emph{The complexity of graph pebbling},
  arxiv.org/abs/math/0503698.

\bibitem{Moews}
David Moews, \emph{Pebbling graphs}, J. Combin. Theory Ser. B \textbf{55}
  (1992), no.~2, 244--252.

\bibitem{Moews_optimal}
\bysame, \emph{Optimally pebbling hypercubes and powers}, Discrete Math.
  \textbf{190} (1998), no.~1-3, 271--276.

\bibitem{Pachter}
Lior Pachter, Hunter~S. Snevily, and Bill Voxman, \emph{On pebbling graphs},
  Proceedings of the Twenty-sixth Southeastern International Conference on
  Combinatorics, Graph Theory and Computing (Boca Raton, FL, 1995), vol. 107,
  1995, pp.~65--80.

\end{thebibliography}

\end{document}